\documentclass[reqno,fleqn,11pt]{amsart}
\usepackage{amsmath}
\usepackage{amssymb}
\usepackage{amsfonts}

\usepackage{graphicx}
\usepackage{amsthm}
\usepackage{enumerate}
\usepackage[T1]{fontenc}
\usepackage{color}

\newcommand{\R}{\mathbb{R}}
\newcommand{\C}{\mathbb{C}}

\newcommand{\Z}{\mathbb{Z}}

\newcommand{\ep}{\epsilon}
\newcommand{\K}{K\"{a}hler}

\newcommand{\dist}{\operatorname{dist}}
\newcommand{\Exp}{\operatorname{Exp}}

\newcommand{\rank}{\operatorname{rank}}

\newtheorem{thm}{Theorem}
\newtheorem{prop}{Proposition}
\newtheorem{lem}[prop]{Lemma}
\newtheorem{cor}[prop]{Corollary}

\newtheorem{rmk}[prop]{Remark}

\newcommand{\fdim}{\hspace*{\fill}$\Box$}
\newcommand{\dimostr}{{\bf Proof. }}

\newcommand{\set}[2]{ \left\{\,#1\,;\,#2\,\right\} }
\newcommand{\gs}{\mathfrak{s}} %%german s
\newcommand{\dirsum}{\sideset{\ }{^{\oplus}}\sum}
\newcommand{\rk}{r}

\begin{document}

\title[K\"{a}hler immersions of homogeneous K\"{a}hler manifolds into...]{K\"{a}hler immersions of
homogeneous K\"{a}hler manifolds into  complex space forms}
\author[A. J. Di Scala, H. Ishi, A. Loi]{Antonio Jose Di Scala, Hideyuki Ishi, Andrea Loi}
\address{Politecnico di Torino, Nagoya University,  Dipartimento di Matematica e Informatica, Universit\`{a} di Cagliari,
Via Ospedale 72, 09124 Cagliari, Italy}
\email{antonio.discala@polito.it; hideyuki@math.nagoya-u.ac.jp; loi@unica.it}
\thanks{Research partially supported by GNSAGA (INdAM), KAKENHI (JSPS) and MIUR (PRIN07, Differential Geometry and Global Analysis)}
\date{September 21, 2010}

\begin{abstract}
In this paper we study the homogeneous \K\ manifolds (h.K.m.) which can be \K\ immersed into  finite or infinite dimensional complex space forms.
On one hand  we completely classify the h.K.m. which  can be \K\ immersed into
a finite or infinite dimensional complex Euclidean or hyperbolic space. Moreover,  we extend  known results about \K\ immersions into  the finite
dimensional complex projective  space to the infinite dimensional setting.
\end{abstract}

\maketitle

\section{Introduction and statements of the main results}
In  this paper we address the following
problem:  {\em classify all
homogeneous \K\ manifolds (\emph{h.K.m.} for short) which admit a   \K\  immersion
into a given finite or   infinite dimensional complex space form}.

\vskip 0.3cm

A \K\ immersion $f:(M, g)\rightarrow (S, g_S)$
from  a \K\ manifold   $(M , g)$  into a complex space form $(S, g_S)$
is a holomorphic map such that $f^*g_S=g$
(here $g$ and $g_S$  denote the  \K\ metrics on $M$ and $S$ respectively).

\vskip 0.3cm

Recall that there are three types,
up to homotheties,
of  complex space forms $(S, g_S)$
according to the sign of
their constant holomorphic sectional curvature:

\begin{itemize}
\item
the complex Euclidean space
$\C^N$, $N\leq\infty$,
with the flat metric denoted by $g_{0}$.
Here ${\C}^{\infty}$ is
the complex  Hilbert space
$\ell^2(\C)$
consisting of
sequences $z_j, j=1\dots, z_j\in {\C}$
such that
$\sum_{j=1}^{+\infty}|z_j|^2<+\infty$.
\item
the complex hyperbolic space
${\C}H^N$,  $N\leq\infty$,
namely the unit ball in $\C ^N$ ($\sum_{j=1}^{N}|z_j|^2<1$)
endowed with the hyperbolic metric
$g_{hyp}$ of
holomorphic sectional curvature being $-4$, whose associated \K\
form  $\omega_{hyp}$ is given by:
\begin{equation}\label{omegahyp}
%\omega_{hyp}=-\frac{i}{2}\partial\bar\partial\log \sum_{j=1}^{N}
%\log (1-|z_j|^2). [corrected]
\omega_{hyp}=-\frac{i}{2}\partial\bar\partial\log(1- \sum_{j=1}^{N}|z_j|^2).
\end{equation}

\item
the complex projective space
${\C}P^N, N\leq\infty$,
with the Fubini--Study metric $g_{FS}$
of  holomorphic sectional curvature being $4$.
If $\omega_{FS}$ denotes the \K\ form associated to $g_{FS}$
then,  in homogeneous
coordinates $[Z_0,\dots, Z_{N}]$,
$\omega_{FS}=\frac{i}{2}\partial\bar\partial\log \sum_{j=0}^{N}
|Z_j|^2$.

\end{itemize}

\vskip 0.3cm

\noindent
{\em Notation.} When we  speak about the \K\ manifold   $\C^N$ (resp. $\C H^N$ or $\C P^N$)
without  mentioning the \K\ metric  we will  always mean  $\C^N$ (resp. $\C H^N$ or $\C P^N$) equipped
with the metric $g_0$ (resp. $g_{hyp}$, $g_{FS}$).

\vskip 0.3cm

Note that, once that a \K\ immersion into a complex space form $(S, g_S)$ is given,
then all other \K\ immersions can be obtained by composing it with  a unitary transformation
of  $(S, g_S)$. This is due to the following celebrated  rigidity theorem due to E. Calabi \cite{Ca53}  which will be of constant use throughout this paper.

\vskip 0.3cm

\noindent
{\bf Theorem (Calabi's rigidity  theorem)}
{\em Let  $f:(M, g)\rightarrow (S, g_S)$  and
$\tilde f:(M, g)\rightarrow (S, g_S)$  be two \K\	immersions into the same  complex space form $(S,
g_S)$.
Then  there exists a unitary transformation $U$ of $(S, g_S)$ such that $f=U\circ\tilde f$.}

\vskip 0.3cm

\subsection{Immersions in $\C^N$ and $\C H^N$}
In the following two theorems   we give a complete  solution  of our problem
when the ambient space is $\C^N$ or $\C H^N$, $N\leq\infty$.
%These theorems   generalize (a), (b), (c) and (d) of Theorem \ref{calforms}  in the Appendix below
%where one can also find the explicit description of the \K\ immersions involved.
In order to state our result
note that  the map $f_n:{\C}H^n\rightarrow l^2({\C})$
given by:
\begin{equation}\label{calhyp}
z=(z_1,\dots z_n)\stackrel{f_n}{\mapsto}
(\dots, \sqrt{\frac{(|j|-1)!}{j!}}z_1^{j_1}\cdots z_n^{j_{n}}, \dots )
\end{equation}
is a \K\ immersion of  ${\C}H^n$ into $l^2({\C})$, i.e. $f_n^*g_0=g_{hyp}$,
(see \cite{Ca53}), where
$|j|=j_1+\dots+j_n$
and $j!=j_1! \cdots j_{n}!$.

%$$\left\{\dots, \frac{\sqrt{( m +j-1)!}}{\sqrt{\pi^n}\sqrt{j_1!\cdots j_{n}!( m -n-1)!}}z_0^{j_0}\cdots z_n^{j_{n-1}},\dots\right\}.$$
%where
%Here we are using the following
%convention:
%we arrange every $n$-tuple of nonnegative
%integers as the sequence
%$m_{j}=(m_{1, j}, m_{2, j},\dots , m_{n, j})
%_{j=0, 1,\dots}$
%such that $m_0=(0, \dots , 0)$,
%$|m_{j}|\leq |m_{j+1}|$, with
%$|m_{j}|=\sum_{\alpha =1}^n m_{\alpha , j}$
%and
%$z^{m_{j}}=\prod_{\alpha =1}^{n}
%(z_{\alpha})^{m_{\alpha , j}}$.
%Further, we order all the $m_j$'s with the same
%$|m_j|$ using the lexicographic order in the variables
%$(z_1, \dots z_n)$.

\begin{thm}\label{Flat-Case}
Let  $(M, g)$ be a $n$-dimensional   h.K.m..
\begin{itemize}
\item [(a)]
If $(M, g)$  can be \K\ immersed into $\C ^N$, $N<\infty$,  then
$(M, g)=\C^n$;
\item [(b)]
if $(M, g)$  can be \K\ immersed into $\ell^2(\C)$,  then
$(M,  g)$  equals
$$\C ^k
\times{\C}H^{n_1}_{\lambda_1}\times\cdots \times {\C}H^{n_\rk}_{\lambda_\rk},$$
where $k+n_1+\cdots +n_\rk=n$, $\lambda_j$, $j=1,\dots , \rk$
are positive real numbers and $\C H^{n_j}_{\lambda_j}=(\C H^{n_j}, \lambda_jg_{hyp}),\ j=1,\dots , \rk$ (hence $\C H^n_{1}=\C H^n$).
\end{itemize}
Moreover, in case (a) (resp. case (b)) the immersion is given, up to a unitary transformation of $\C ^N$ (resp. $\ell^2 (\C))$,  by the linear inclusion $\C^n\hookrightarrow \C ^N$
(resp. by
$(f_0, f_1,\dots ,f_\rk)$,
where $f_0$ the  linear inclusion
$\C ^k\hookrightarrow \ell^2(\C)$
and   each $f_j:{\C}H^{n_j}\rightarrow \ell^2({\C})$
is
$\sqrt{\lambda_j}$ times
the map (\ref{calhyp})).
\end{thm}

\begin{thm}\label{Negative-case}
Let  $(M, g)$  be a $n$-dimensional  h.K.m..
Then if $(M, g)$  can be \K\ immersed into $\C H^N$, $N\leq\infty$,  then
$(M, g)=\C H^n$
and the immersion is given,
up to  a unitary transformation of $\C H^N$,
by the linear inclusion $\C H^n\hookrightarrow \C H^N$
\end{thm}

\begin{rmk}\rm
Since a \K\ immersion is minimal, an alternative proof of  (1) in  Theorem \ref{Flat-Case} when $N<\infty$
follows by the work of
A. J. Di Scala \cite{DS02}.
%which asserts that  if  a $n$-dimensional h.K.m.  admits a minimal  immersion into
%$\C ^N$  then $M=\C ^n$.
\end{rmk}

\begin{rmk}\rm\label{assertion2}
Assertion (2) in Theorem  \ref{Flat-Case}
is a generalization  to arbitrary h.K.m.  of  Theorem 3.3  in \cite{DL07}
where the first and the third authors  proved that  a  bounded
symmetric domain which can be \K\ immersed into $\ell ^2 (\C)$ has necessarily  rank  one.
Actually, the method of the present paper,  when applied to bounded symmetric domains, provides
us with  an alternative and more elegant   proof of
this result  (cfr. Remark \ref{simpleproof} below).
\end{rmk}

\vskip 0.3cm

\subsection{Immersion in $\C P^N$}
There exists a large class (cfr. Conjecture 1 below)  of  h.K.m. which can be \K\ immersed into $\C P^{N}$.
In this paper  a  K\"ahler metric $g$ on a complex  manifold $M$ will be called  {\em projectively induced} if there exists an immersion  $f:M \rightarrow \mathbb{C}P^{N}$, $N\leq\infty$,  such that $f^*g_{FS}=g$.
An obvious   necessary condition for $g$ to be projectively induced is that its associated \K\ form $\omega$ is integral
i.e.  it   represents the first Chern class $c_{1}(L)$ in $H^2(M, \Z)$
of a holomorphic line bundle $L\rightarrow M$. Indeed $L$ can be taken as the pull-back of the hyperplane line  bundle on $\C P^N$
whose first Chern class is given by $\omega_{FS}$.
Notice that if $\omega$ is an exact form (e.g. when  $M$ is contractible) then $\omega$ is obviously
integral  since its  second cohomology class vanishes.

Other (less obvious) conditions are expressed by the following theorem and its corollary which represent our first result
about  projectively induced \K\ metrics.

\vskip 0.3cm

\begin{thm}\label{necessary}
Assume that a h.K.m. $(M, g)$ admits a K\"ahler immersion $f : M \rightarrow \C P^{N}$, $N\leq\infty$. Then
$M$ is simply-connected and $f$ is injective.
\end{thm}

\begin{cor}\label{corolnecessary} Let $(M, g)$ be a complete and  locally h.K.m..
Assume that  $f: (M, g)\rightarrow \C P^N$, $N\leq\infty$,  is a \K\ immersion.
Then $(M, g)$ is a h.K.m..
\end{cor}

When the dimension of the  ambient space is finite, i.e. $(S, g_S)=\C P^N$, $N<\infty$,
$M$ is forced to be  compact and a proof of Theorem \ref{necessary} is
well-known  by the work of M. Takeuchi \cite{TA78}.
In this case he also provides a complete classification of all compact
h.K.m. which can be \K\ immersed into $\C P^N$
by making use of the representation theory of semisimple Lie groups.
Viceversa,  it is not hard to see   that
if a {\em compact}  \K\ manifold can be \K\ immersed into $\C P^{\infty}$
then  it can also be \K\ immersed into $\C P^N$ with $N<\infty$.
%Indeed, note  first  that   the immersion $f:M\rightarrow \C P^{\infty}$  is given by global holomorphic sections of the line bundle $L$ obtained as the pull-back
%of the hyperplane bundle of $\C P^{N}$. Secondly,  one  can also assume that  the immersion is full, i.e. $f(M)$ is not contained into a proper complex projective space of $\C P^%{\infty}$;  this means that the immersion is built using linearly independent sections of $L$.
%On the other hand,   the dimension of the  space of global holomorphic sections of $L$
%is finite since $M$ is compact and this proves the claim.

We believe that, up to homotheties,   {\em any} simply-connected  h.K.m. such that its associated  \K\  form is integral can be \K\ immersed into $\C P^N$,
with $N\leq\infty$. This is  expressed by the
following conjecture.

\vskip 0.3cm

\noindent
{\bf Conjecture 1:}  {\em
Let $(M, g)$ be a simply-connected h.K.m.  such that its  associated \K\ form  $\omega$ is integral.
Then  there exists $\lambda_0 \in \mathbb{R}^+$ such that $\lambda_0 g$  is  projectively induced.
}

\vskip 0.3cm

The integrality of  $\omega$  in the conjecture  is important since there exist  simply-connected h.K.m. $(M, \omega)$ such that $\lambda \omega$ is not integral for any $\lambda \in \mathbb{R}^+$ (take,  for example, $(M, g) = (\mathbb{C}P^1, g_{FS}) \times (\mathbb{C}P^1, \sqrt{2}g_{FS})$).
%For a non simply-connected example take $M = (\mathbb{C}/\Gamma, \omega_{0}) \times (\mathbb{C}P^1, \sqrt{2}\omega_{FS})$.
Observe also that  there exist simply-connected (even contractible) h.K.m. $(M, g)$ such that $\omega$ is an integral form but
$g$ is not projectively induced.
In order to describe such an example we  recall the following result (see Theorem 2 in  \cite{LZ09}).

\vskip 0.3cm

\noindent
{\bf Theorem A.}
Let $g_B$ be the Bergman metric  of an irreducible
Hermitian symmetric space of noncompact type $\Omega$.
Then  $\lambda g_B$
is projectively induced if and only $\lambda\gamma$ belongs to $W(\Omega)\setminus \{0\}$, where $\gamma$ denotes the genus of $\Omega$ and $W(\Omega)$ its Wallach set.

\vskip 0.3cm

It turns out (see Corollary $4.4$ p. 27 in  \cite{AR95} and references therein)
that $W(\Omega)$ consists only of real numbers and depends   on two of the domain's invariants, denoted  by  $a$ (strictly positive natural  number) and  $r$ (the rank of $\Omega$).
More precisely we have
\begin{equation}\label{wallachset}
W(\Omega)=\left\{0,\,\frac{a}{2},\,2\frac{a}{2},\,\dots,\,(r-1)\frac{a}{2}\right\}\cup \left((r-1)\frac{a}{2},\,\infty\right).
\end{equation}
The set $W_d=\left\{0,\,\frac{a}{2},\,2\frac{a}{2},\,\dots,\,(r-1)\frac{a}{2}\right\}$ and the interval $W_c= \left((r-1)\frac{a}{2},\,\infty\right)$
are called respectively  the {\em discrete} and {\em continuous} part   of the Wallach set of the domain
$\Omega$.
Observe that when $r=1$, namely $\Omega$ is the complex hyperbolic space $\C H^n$,
then $g_B=(n+1)g_{hyp}$.
In this case (and only in this case)   $W_d=\{0\}$ and $W_c=(0, \infty)$.
If  $\rank (\Omega)=r\geq 2$
and $0<\lambda < \frac{a}{2\gamma}$
it follows by Theorem A that
$\lambda g_B$ is not projectively induced
and its associated \K\ form $\lambda \omega_B$ is integral (since $\Omega$ is contractible).
This provides us with the  desired example.

\vskip 0.3cm
Notice also that from Theorem A it follows that  the only irreducible bounded symmetric domain where  $\lambda g_B$
is projectively induced for all $\lambda >0$ is the complex hyperbolic space.
In the following theorem,  which represents our last result,  we generalize  this fact to
any homogeneous bounded domain (h.b.d. for short).
This will be a key ingredient in the proof of Theorem \ref{Flat-Case}.

\begin{thm}\label{thmsmall}
Let $(\Omega, g)$ be a $n$-dimensional  h.b.d..
%Assume that there exists a positive real number  $\lambda_0>0$
%such that $\lambda g$ is  projectively induced for all real numbers
%$0<\lambda \leq \lambda_0$. [corrected]
The metric $\lambda g$ is  projectively induced for all
$\lambda >0$
if and only if
\begin{equation} \label{eqn:direct_prod}
(\Omega, g) ={\C}H^{n_1}_{\lambda_1}\times\cdots \times {\C}H^{n_\rk}_{\lambda_\rk},
\end{equation}
where $n_1+\cdots +n_r=n$,  $\lambda_j$, $j=1,\dots , \rk$
are positive real numbers and
$\C H^{n_j}_{\lambda_j}=(\C H^{n_j}, \lambda_jg_{hyp}),\ j=1,\dots , \rk$.
\end{thm}

\vskip 0.3cm
\noindent
The paper contains another section dedicated to the proofs of  our main results.

\vskip 0.3cm

\noindent
{\bf Aknowledgments:}
The second and third author would like to thank {\em Politecnico of Torino} for the wonderful hospitality
in their research stays in January 2010.

\section{Proof of the main results}

The  basic ingredient for the proof of  our results
is the following solution due to
J. Dorfmeister and   K. Nakajima \cite{DN88} of the fundamental conjecture on h.K.m..

\vskip 0.3cm

\noindent
{\bf Theorem FC}
{\em A h.K.m. $(M, g)$  is the total space of a holomorphic
fiber bundle over a h.b.d. $\Omega$ in which the fiber ${\mathcal F} ={\mathcal E} \times {\mathcal C}$
is (with the
induced K\"ahler metric) the  \K\ product of a flat homogeneous K\"ahler manifold ${\mathcal E}$ and a
compact simply-connected homogeneous K\"ahler manifold ${\mathcal C}$.}

\vskip 0.3cm

In order to prove Theorem \ref{Flat-Case}
recall that complete connected totally geodesic submanifolds of $\mathbb{R}^n$ are affine subspaces $p + \mathbb{W}$,
where $p \in \mathbb{R}^n$ and $\mathbb{W} \subset \mathbb{R}^n$ is a vector subspace.
We need the following result  from \cite{AD03} which we include here for completeness.

\begin{lem} \label{parallel} Let $G$ be a connected  Lie subgroup of
isometries of the Euclidean space $\mathbb{R}^n$. Let $G.p = p + \mathbb{V}$ and $G.q = q + \mathbb{W}$ be two totally geodesic $G$-orbits.
Then $\mathbb{V} = \mathbb{W}$, i.e. $G.p$ and $G.q$ are parallel affine subspaces of $\mathbb{R}^n$.
\end{lem}
\noindent
\dimostr
We can assume that $p = 0 \in \mathbb{R}^n$ and that $p,q$ are the points that realize the distance between both orbits
$G\cdot p$, $G\cdot q$, i.e.
$\dist(p,q) = \dist(G\cdot p,G\cdot q)$. Let $\gamma(t)=tq$ be the geodesic
that realizes the distance between $q$ and $\mathbb{V}$. So the vector $q$
is perpendicular to any $G$-orbit $G_t=G \cdot \gamma(t)$ $t \in
{\mathbb R}$. Let $X = x^*$ be any Killing vector field of $G$ and $\Exp(tX)$
its associated one-parameter group of isometries. Define
$h:I\times \mathbb{R }\rightarrow \mathbb{R}^n$ by $h_s(t):=\Exp(sX) \cdot
\gamma(t)$. Note that $X ( h_s(t)) = \frac{\partial h}{\partial
s}$ and that, for a fixed $s$, $h_s(t)$ is a geodesic.

Let $A_t$ be the shape operator at the point $\gamma(t)$ of the
orbit $G \cdot \gamma(t)$ in the direction of $\dot \gamma (t)$.
Define $f(t) := - \langle A_t(X ( \gamma(t))) ,X (\gamma(t))
\rangle =
\langle \frac{D}{\partial s}\frac{\partial h}{\partial t},
X(h_s(t)) \rangle \mid_{s=0}$. We have
$$\frac{d}{dt}f(t) = \langle \frac{D}{\partial t}
\frac{D}{\partial s}\frac{\partial h}{\partial t}, X(h_s(t))
\rangle \mid_{s=0} + \langle \frac{D}{\partial s} \frac{\partial
h}{\partial t}, \frac{D}{\partial t}X(h_s(t)) \rangle \mid_{s=0}
$$
$$  = \langle \frac{D}{\partial s}
\frac{D}{\partial t}\frac{\partial h}{\partial t}, X(h_s(t))
\rangle \mid_{s=0}
+ \langle \frac{D}{\partial t}\frac{\partial h}{\partial s},
\frac{D}{\partial t}X (h_s(t)) \rangle \mid_{s=0} $$
$$ =  \| \nabla_{\dot{\gamma}(t)}
(X(\gamma(t))) \| ^2 .$$

Since $f(0)=0$ because $G \cdot p$ is totally geodesic, we get
$$f(1)= - \langle A_t(X( q)), X( q) \rangle \geq 0 \, .$$ Hence $A_1$ is
negative definite and since $G \cdot q$ is totally geodesic, any Killing
vector field $X$ is parallel along $\gamma(t)$. We can write
$\Exp(sX) \cdot p = e^{s\overline{X}}(p-c)+c+sd$, where $\overline{X}$ is the
projection of $X$ into ${\mathfrak{so}}_n$, $d \in \ker (\overline{X})$ and
$c \in \ker (\overline{X})^{\perp}$. Then a Killing vector field $X$ is
parallel along $\gamma(t)$ if and only if $q \in \ker (\overline{X})$.
Thus $\Exp(sX) \cdot q = q + \Exp(sX) \cdot p $ which implies that
$$\mathbb{V} = T_p (G\cdot p) \subset T_q (G\cdot q) = \mathbb{W} \, .$$
Reversing the role of $ \mathbb{V}$ and  $\mathbb{W}$ the same argument  yields $\mathbb{W} \subset \mathbb{V}$.
This  completes the proof of the lemma. $\Box$

\vskip 0.3cm

\noindent
{\bf Proof of Theorem \ref{Flat-Case}.}
\noindent
Assume that there exists a K\"ahler immersion $f: M \rightarrow \C^N$. By Theorem FC
and by the fact that a h.b.d. is contractible
we get that  $M = \mathbb{C}^k\times \Omega$ as a complex manifold since,
by the maximum principle, the fiber ${\mathcal F}$ cannot contain a compact manifold. Let $M = G/K$ be the homogeneous realization  of $M$ (so the metric
$g$ is $G$-invariant).
It follows  again by Theorem FC that there exists $L \subset G$ such that the $L$-orbits are the fibers of the fibration $\pi: M = G/K \rightarrow \Omega = G/L $. Let $F_p, F_q$ be the fibers over $p,q \in \Omega$. We claim that $f(F_p)$ and $f(F_q)$ are parallel affine subspaces of $\C^N$.
Indeed, by Calabi's rigidity  $f(F_p)$ and $f(F_q)$ are affine subspaces of $\C^N$ since both $F_p$ and $F_q$ are flat \K\ manifolds  of $\mathbb{C}^n$.
Moreover, Calabi rigidity theorem implies the existence of a  morphism of groups $\rho: G \rightarrow Iso_{\mathbb{C}}(\C^N) = \mathrm{U}(\C^N) \ltimes \C^N$ such that
$f(g\cdot x) = \rho(g)f(x)$
for all $g \in G, x \in M$.
Let $W_{p,q}$ be the affine subspace generated by $f(F_p)$ and $f(F_q)$. Since both $f(F_p)$ and $f(F_q)$ are $\rho(L)$-invariant it follows that $W_{p,q}$ is also $\rho(L)$-invariant. Indeed, for any $g \in L$ the isometry $\rho(g)$ is an affine map and so must preserve the affine space generated by $f(F_p)$ and $f(F_q)$. Observe that $W_{p,q}$ is a finite dimensional complex Euclidean space,  $\rho(L)$ acts on $W_{p,q}$
and $f(F_p)$ and $f(F_q)$ are  two complex totally geodesic orbits in $W_{p,q}$. Then,  by  Lemma \ref{parallel}, we get that $f(F_p)$ and $f(F_q)$ are parallel affine subspaces of $W_{p,q}$ and hence  of  $\C^N$. Since $p,q \in \Omega$ are two arbitrary points it follows that $f(M)$ is a K\"ahler product. Thus $M = \mathbb{C}^{k}\times \Omega$ is a K\"ahler product of  homogeneous K\"ahler manifolds. Using again the fact  $M$ can be \K\ immersed into $\C^N$ it follows that the h.b.d.  $\Omega$   can be \K\ immersed into $\C^N$.
If one denotes by $\varphi$ this immersion  and by
$g_{\Omega}$ the homogeneous
\K\ metric of $\Omega$,
it follows that  the map $\sqrt{\lambda}\varphi$ is a \K\ immersion of $(\Omega, \lambda g_{\Omega})$ into $\C^N$. Therefore, by Theorem  14 in \cite{Bo47}, $\lambda g_{\Omega}$ is projectively  induced for all $\lambda >0$ and
Theorem \ref{thmsmall} yields
$$(M, g) = \C ^k
\times {\C}H^{n_1}_{\lambda_1}\times\cdots \times {\C}H^{n_\rk}_{\lambda_\rk},$$
where $k+n_1+\cdots +n_\rk=n$ and $\lambda_j$, $j=1,\dots , \rk$
are positive real numbers.
If the dimension $N$ of the ambient space $\C^N$ is finite then
$M =\C^n$ since there cannot exist a \K\ immersion of
$(\C H^{n_j}, \lambda_j g_{hyp})$
into $\C^N$, $N<\infty$ (see \cite{Ca53})
and this proves (a).
The last part of Theorem \ref{Flat-Case}
is a consequence of   Calabi's rigidity theorem together with
Lemma 3.1 in \cite{DL07}
which asserts that
a \K\ map $f: M \times M' \rightarrow \C ^N$, $N\leq\infty$,
from a product $M \times M'$ of two \K\
manifolds is a product, i.e. $f(p,q) = (f_1(p),f_2(q))$ where
$f_1:M \rightarrow \C ^N$  and $f_2:M' \rightarrow \C ^N$ are
\K\ maps.
$\Box$

\begin{rmk}\label{simpleproof}\rm
As we have already pointed, Theorem \ref{thmsmall}, which is an important step in the
proof of the Theorem \ref{Flat-Case},
is  a straightforward  consequence of Theorem A above when the h.K.m. is  a  bounded symmetric domain.
Therefore the last part of Theorem \ref{Flat-Case} provides  an alternative proof of  Theorem 3.3  in \cite{DL07}
without the use of Calabi's diastasis function (cfr. Remark \ref{assertion2}).
\end{rmk}

\vskip 0.3cm

In order to prove Theorem \ref{Negative-case} we need the following lemma.

\begin{lem}\label{lemmazedda1}
If  a \K\ manifold $(M, g)$    can be  K\"ahler immersed into $\C H^N$, $N\leq\infty$,  then it can also be
K\"ahler immersed into $\ell^2(\mathbb{C})$.
\end{lem}
\noindent
\dimostr
Let $f$ be the  \K\ immersion of $(M, g)$ into $\C H^N$. If $N<\infty$  then the map  $f_n\circ f : (M, g)\rightarrow \ell ^2(\C)$,
where $f_n$ is given by (\ref{calhyp}), is a  \K\ immersion.  If $N=\infty$,
it follows by (\ref{omegahyp}) in the introduction  that $\Phi =-\log (1-\sum_{j=1}^{\infty}|\phi_j|^2)
= \sum_{k=1}^{\infty}(\sum_{j=1}^{\infty}|\phi_j|^2)^k$ is a \K\ potential
for the metric $g$, i.e. $\frac{i}{2}\partial\bar\partial \Phi =\omega$, where $\omega$ is the \K\ form associated to the metric $g$
and  the $\phi_j$'s are the components of $f$.
Hence $\Phi =\sum_{j=1}^{\infty}|h_j|^2$
for suitable  holomorphic functions $h_j$, $j=1, 2, \dots$ on $M$  and
the map $h= (\dots, h_j, \dots):(M, g)\rightarrow  \ell ^2(\C)$ is the desired \K\ immersion.
\fdim

\vskip 0.3cm

\noindent
{\bf Proof of Theorem \ref{Negative-case}.}
If  a h.K.m. $(M, g)$  can be
\K\ immersed into $\C H^{N}$,
$N\leq\infty$,
then, by Lemma \ref{lemmazedda1}  it can also be
K\"ahler immersed into $\ell^2(\mathbb{C})$.
By Theorem \ref{Flat-Case},
$(M, g)$ is then   a \K\ product of complex space forms, namely
$$(M,  g)=\C ^k
\times {\C}H^{n_1}_{\lambda_1}\times\cdots \times {\C}H^{n_\rk}_{\lambda_\rk}.$$
Then  the conclusion follows by the fact that $\C^k$ cannot be \K\ immersed into
$\C H^N$ for all $N\leq\infty$ (see \cite{Ca53}), by Calabi's rigidity theorem
and by  Theorem 2.11 in \cite{AD03}
which  shows that there are not \K\ maps
from a product $M \times M'$ of \K\ manifolds into ${\C}
H^{N}$, $N\leq\infty$, (the proof in \cite{AD03} is given for $N<\infty$  but it
extends without any substantial change
to the infinite dimensional case).
$\Box$

\vskip 0.3cm

\noindent
{\bf Proof of Theorem \ref{necessary}.}
Theorem FC and the fact that a h.b.d. is contractible imply that $M$ is a {\em complex}
product $\Omega \times {\mathcal F}$, where  ${\mathcal F}={\mathcal E}\times {\mathcal C}$ is a \K\ product  of a flat  \K\ manifold ${\mathcal E}$ \K\ embedded into $(M, g)$ and a simply-connected h.K.m. ${\mathcal C}$.
We claim that ${\mathcal E}$ is simply-connected and hence   $M=\Omega\times {\mathcal E}\times {\mathcal C}$ is simply-connected.
In order to prove our claim notice that
${\mathcal E}$ is the \K\ product  $\C^k\times T_1\times\cdots\times T_s$, where $T_j$ are flat complex tori.
So one needs to show that each $T_j$  reduces to a point.
If, by a contradiction,  the dimension of one of this tori, say $T_{j_0}$  is not zero, then by composing
the \K\  immersion of  $T_{j_0}$ in $(M, g)$
with the immersion  $f:M\rightarrow \C P^N$
we would get  a \K\ immersion of
$T_{j_0}$ into $\C P^N$ in contrast with a well-known result of Calabi \cite{Ca53} (see also  Lemma 2.2 in \cite{TA78}).
In order to prove that $f$ is injective we first observe that, by  Calabi's rigidity theorem,  $f(M)$ is still a h.K.m..
Then, by  the first part of the theorem,  $f(M)\subset\C P^{N}$ is simply-connected. Moreover,  since $M$ is complete and  $f: M \rightarrow f(M)$ is a local isometry,  it is a covering map (see,  e.g.,  Lemma 3.3 p. 150 in \cite{DC92}) and  hence injective.$\Box$

\vskip 0.3cm

\noindent
{\bf Proof of Corollary  \ref{corolnecessary}.}
Let $\pi:\tilde M\rightarrow M$ be the universal  covering map.
Then $(\tilde M, \tilde g)$ is a  h.K.m. and, by Theorem \ref{necessary},
$f\circ\pi :\tilde M\rightarrow \C P^n$ is injective.
Therefore $\pi$ is injective,
and since it is a covering map,
it defines a holomorphic isometry between $(\tilde M, \tilde g)$
and $(M, g)$.
$\Box$

\vskip 0.3cm

\noindent
{\bf Proof of Theorem \ref{thmsmall}.}
First we find a global potential of the homogeneous K\"ahler metric $g$
on the domain $\Omega$ following Dorfmeister \cite{D85}.
By \cite[Theorem 2 (c)]{D85},
there exists a split solvable Lie subgroup
$S \subset \mathrm{Aut}(\Omega, g)$
acting simply transitively on the domain $\Omega$.
Taking a reference point $z_0 \in \Omega$,
we have a diffeomorphism
$S \owns s \overset{\sim}{\mapsto} s \cdot z_0 \in \Omega$,
and
by the differentiation, we get the linear isomorphism
$\gs := \mathrm{Lie}(S) \owns X
\overset{\sim}{\mapsto} X \cdot z_0 \in T_{z_0}\Omega \equiv \C^n$.
Then the evaluation of the K\"ahler form $\omega$ on $T_{z_o}\Omega$
is given by
$\omega(X\cdot z_o, Y \cdot z_0) = \beta([X,Y])\,\,
(X, Y \in \gs)$
with a certain linear form $\beta \in \gs^*$.
Let $j : \gs \to \gs$ be the linear map defined in such a way that
$(jX) \cdot z_0 = \sqrt{-1} (X \cdot z_0)$ for $X \in \gs$.
We have
$\Re g(X \cdot z_0,\,Y \cdot z_0) = \beta([jX, Y])$
for $X, Y \in \gs$,
and the right-hand side defines a positive inner product on $\gs$.
Let $\mathfrak{a}$ be the orthogonal complement of $[\gs, \gs]$ in $\gs$
with respect to the inner product.
Then $\mathfrak{a}$ is a commutative Cartan subalgebra of $\gs$.
Define $\gamma  \in \mathfrak{a}^*$ by
$\gamma(C) := -4 \beta (jC)\,\,\,(C \in \mathfrak{a})$,
and we extended $\gamma$ to $\gs = \mathfrak{a} \oplus [\gs, \gs]$ by the zero-extension.
Keeping the diffeomorphism between $S$ and $\Omega$ in mind,
we define a positive smooth function $\Psi$ on $\Omega$ by
$$ \Psi((\exp X) \cdot z_0) = e^{-\gamma(X)} \,\,\, (X \in \gs). $$
From the argument in \cite[pp. 302--304]{D85},
we see that
\begin{equation} \label{eqn:globalpotential}
\omega = \frac{i}{2}\partial\bar\partial\log \Psi.
\end{equation}
It is known that
there exists a unique kernel function $\tilde{\Psi} : \Omega \times \Omega \to \C$
such that (1) $\tilde{\Psi}(z,z) = \Psi(z)$ for $z \in \Omega$
and (2) $\tilde{\Psi}(z,w)$ is holomorphic in $z$ and anti-holomorphic in $w$
(cf. \cite[Proposition 4.6]{I99}).
Let us observe that the metric $g$ is projectively induced
if and only if
$\tilde{\Psi}$ is a reproducing kernel of a Hilbert space
of holomorphic functions on $\Omega$.
Indeed,
if $f : \Omega \to \C P^N\,\,(N \le \infty)$ is a K\"ahler immersion
with
$f(z) = [\psi_0(z) : \psi_1(z) : \cdots ]\,\,(z \in \Omega)$
its homogeneous coordinate expression,
then we have
$\omega = \frac{i}{2} \partial\bar\partial\log \sum_{j=0}^{N}|\psi_j|^2$.
Comparing (\ref{eqn:globalpotential}) with it,
we see that
there exists a holomorphic function $\phi$ on $\Omega$ for which
$\Psi = |e^{\phi}|^2 \sum_{j=0}^N |\psi_j|^2$.
By analytic continuation,
we obtain
$\tilde{\Psi}(z,w) = e^{\phi(z)} \overline{e^{\phi(w)}} \sum_{j=0}^N \psi_j(z) \overline{\psi_j(w)}$
for $z, w \in \Omega$.
For any $z_1, \dots, z_m \in \Omega$ and $c_1, \dots, c_m \in \C$,
we have
\begin{align*}
\sum_{p,q=1}^m c_p \bar{c}_q\tilde{\Psi}(z_p,z_q)
&= \sum_{p,q=1}^m c_p \bar{c}_q e^{\phi(z_p)} \overline{e^{\phi(z_q)}} \sum_{j=0}^N \psi_j(z_p) \overline{\psi_j(z_q)} \\
&= \sum_{j=0}^N |\sum_{p=1}^m c_p e^{\phi(z_p)}\psi_j(z_p)|^2 \ge 0.
\end{align*}
Thus the matrix $(\tilde{\Psi}(z_p,z_q))_{p,q} \in \mathrm{Mat}(m,\C)$ is always
a positive Hermitian matrix.
Therefore $\tilde{\Psi}$ is a reproducing kernel of a Hilbert space
(see \cite[p. 344]{Ar50}).

On the other hand,
if $\tilde{\Psi}$ is a reproducing kernel of a Hilbert space
$\mathcal{H} \subset \mathcal{O}(\Omega)$,
then by taking an orthonormal basis $\{\psi_j\}_{j=0}^N$ of $\mathcal{H}$,
we have a K\"ahler immersion $f : M \owns z \mapsto [\psi_0(z) : \psi_1(z) : \cdots] \in \C P^N$
because we have $\Psi(z) = \tilde{\Psi}(z,z) = \sum_{j=0}^N |\psi_j(z)|^2$.
Note that there exists no point $a \in \Omega$ such that $\psi_j(a) = 0$ for all $1 \le j \le N$
since $\Psi(z) = \sum_{j=0}^N |\psi_j(z)|^2$ is always positive.

The condition for $\tilde{\Psi}$ to be a reproducing kernel is described in \cite{I99}.
In order to apply the results,
we need a fine description of the Lie algebra $\gs$ with $j$
due to Piatetskii-Shapiro \cite{PS69}.
Indeed, it is shown in \cite[Chapter 2]{PS69} that
the correspondence between the h.b.d. $\Omega$ and the structure of $(\gs, j)$
is one-to-one up to natural equivalence.
For a linear form $\alpha$ on the Cartan algebra $\mathfrak{a}$,
we denote by $\gs_{\alpha}$ the root subspace
$\set{X \in \gs}{[C,X] = \alpha(C)X \,\, (\forall C \in \mathfrak{a})}$
of $\gs$.
The number $r := \dim \mathfrak{a}$
is nothing but the rank of $\Omega$.
Thanks to \cite[Chapter 2, Section 3]{PS69},
there exists a basis
$\{\alpha_1, \dots, \alpha_r\}$
of $\mathfrak{a}^*$
such that
$\gs = \gs(0) \oplus \gs(1/2) \oplus \gs(1)$ with
\begin{align*}
\gs(0) &= \mathfrak{a} \oplus
          \dirsum_{1 \le k < l \le r} \gs_{(\alpha_l - \alpha_k)/2}, \quad
\gs(1/2) = \dirsum_{1 \le k \le r} \gs_{\alpha_k /2}, \\
\gs(1) &= \dirsum_{1 \le k \le r} \gs_{\alpha_k}
        \oplus
        \dirsum_{1 \le k < l \le r} \gs_{(\alpha_l + \alpha_k)/2}.
\end{align*}
If $\{A_1, \dots, A_r\}$ is the basis of $\mathfrak{a}$ dual to $\{\alpha_1, \dots, \alpha_r\}$,
then $\gs_{\alpha_k} = \R jA_k$.
Thus $\gs_{\alpha_k}\,\,(k=1, \dots, r)$ is always one dimensional,
whereas other root spaces $\gs_{\alpha_k/2}$ and $\gs_{(\alpha_l \pm \alpha_k)/2}$ may be $\{0\}$.
Since $\{\alpha_1, \dots, \alpha_r\}$ is a basis of $\mathfrak{a}^*$,
the linear form $\gamma \in \mathfrak{a}^*$ is written as
$\gamma = \sum_{k=1}^r \gamma_k \alpha_k$
with unique $\gamma_1, \dots, \gamma_r \in \R$.
Since $j A_k \in \gs_{\alpha_k}$,
we have
\begin{align*}
\gamma_k = \gamma(A_k) = -4 \beta (jA_k) = -4 \beta([A_k, jA_k]) = 4 \beta([jA_k, A_k])
\end{align*}
and the last term equals $4 g (A_k \cdot z_0, A_k \cdot z_0)$.
Thus we get $\gamma_k >0$.

For $\ep = (\ep_1, \dots, \ep_r) \in \{0,1\}^r$,
put
$q_k(\ep) := \sum_{l>k} \ep_l \dim \gs_{(\alpha_l- \alpha_k)/2}
\,\,\,(k=1, \dots, r)$.
Define
$$
\mathfrak{X}(\ep) := \set{(\sigma_1, \dots, \sigma_r) \in \C^r}
{\begin{aligned}
\sigma_k &> q_k(\ep) /2 \quad (\ep_k = 1)\\
\sigma_k &= q_k(\ep) /2 \quad (\ep_k = 0)
\end{aligned}},
$$
and $\mathfrak{X} := \bigsqcup_{\ep \in \{0,1\}^r} \mathfrak{X}(\ep)$.
By \cite[Theorem 4.8]{I99},
$\tilde{\Psi}$ is a reproducing kernel if and only if
$\underline{\gamma} := (\gamma_1, \dots, \gamma_r)$
belongs to $\mathfrak{X}$.
We denote by $W(g)$ the set of $\lambda >0$ for which $\lambda g$ is projectively induced.
Since
the metric $\lambda g$ corresponds to the parameter $\lambda \underline{\gamma}$,
we see that $\lambda g$ is projectively induced if and only if $\lambda \underline{\gamma} \in \mathfrak{X}$.
Namely we obtain
$$
W(g) = \set{\lambda >0}{\lambda \underline{\gamma} \in \mathfrak{X}},
$$
and the right-hand side is considered in \cite{I10}.
Put $q_k = \sum_{l>k} \dim \gs_{(\alpha_l - \alpha_k)/2}$ for $k=1, \dots, r$.
Then \cite[Theorem 15]{I10} tells us that
$$
W(g) \cup \{0\} \subset \set{\frac{q_k}{2\gamma_k}}{k=1, \dots, r} \cup (c_0, +\infty),
$$
where $c_0 := \max \set{\frac{q_k}{2\gamma_k}}{k=1, \dots, r}$.

Now assume that $\lambda g$ is projectively induced for all $\lambda >0$.
Then we have $c_0 = 0$, so that $\dim \gs_{(\alpha_l-\alpha_k)/2} = 0$ for all $1 \le k < l \le r$.
In this case, we see that $\gs$ is a direct sum of ideals
$\gs_k := j \gs_{\alpha_k} \oplus \gs_{\alpha_k/2} \oplus \gs_{\alpha_k}
\,\,\,(k=1, \dots, r),$
which correspond to the hyperbolic spaces $\C H^{n_k}$ with $n_k = 1 + (\dim_{\alpha_k/2})/2$
(\cite[pp. 52--53]{PS69}).
Therefore the Lie algebra $\gs$ corresponds to the direct product
$\C H^{n_1} \times \cdots \times \C H^{n_\rk}$,
which is biholomorphic to $\Omega$ because the homogeneous domain
$\Omega$ also corresponds to $\gs$.
Hence (\ref{eqn:direct_prod}) holds and Theorem 4 is verified.
\fdim

%On the other hand,
% assume that (\ref{eqn:direct_prod}) holds.
%Then we get again $\dim \gs_{(\alpha_l - \alpha_k)/2} = 0$, so that $\mathfrak{X} = [0, +\infty)^r$.
%Therefore for all $\lambda >0$
% we have $\lambda \underline{\gamma} \in \mathfrak{X}$, which means that
% $\lambda g$ is always projectively induced.

\end{document}